\documentclass[entropy,article,submit,moreauthors,pdftex]{Definitions/mdpi}
\preto{\abstractkeywords}{\nolinenumbers}
\firstpage{1}
\pubvolume{xx}
\issuenum{x}
\articlenumber{x}
\pubyear{2022}
\copyrightyear{2022}

\pdfoutput=1

\usepackage{amsmath, amssymb, amsthm, amsfonts, amsxtra}
\usepackage{graphicx,xcolor}
\usepackage{tikz}
\usepackage[all]{xy}

\graphicspath{{Figs/}}

\setitemize{parsep=6pt,itemsep=0pt,leftmargin=*,labelsep=5.5mm}
\setenumerate{parsep=6pt,itemsep=0pt,leftmargin=*,labelsep=5.5mm} 
\setlist[description]{itemsep=0mm}   
\usepackage{environ}
\NewEnviron{myequation}{%
\begin{equation} 
\scalebox{0.95}{$\BODY$}
\end{equation}
}

\newcommand{\lbar}[1]{ \overline{#1} }

\newcommand{\set}[2]{ \left\{\,#1\,;\,#2\,\right\} }
\newcommand{\cpling}[2]{ \langle #1,#2 \rangle}
\newcommand{\transp}[1]{ #1^{\top}}

\newcommand{\nc}{\newcommand}
\nc{\iu}{i}
\nc{\integer}{\mathbb{Z}}
\nc{\real}{\mathbb{R}}
\nc{\complex}{\mathbb{C}}

\nc{\Zsp}{\mathcal{Z}}
\nc{\Lsp}{\mathcal{L}}
\nc{\Ssp}{\mathcal{S}}
\nc{\Pcone}{\mathcal{P}}
\nc{\Vertex}{\mathcal{V}}
\nc{\Graph}{\mathcal{G}}

\Title{Graphical Gaussian models associated to a homogeneous graph with permutation symmetries}

\Author{
{Piotr Graczyk} $^1$, 
{Hideyuki Ishi} $^2$, 
{Bartosz Ko{\l}odziejek} $^3$
} 

\AuthorNames{{Piotr Graczyk}, {Hideyuki Ishi}, {Bartosz Ko{\l}odziejek}}

\address{%
1- Universit\'e d'Angers, France,\quad 
{\it piotr.graczyk@univ-angers.fr}
\\
2- Osaka Metropolitan University, Japan,\quad 
{\it hideyuki-ishi@omu.ac.jp} 
\\
3- Warsaw University of Technology, Poland,\quad 
{\it b.kolodziejek@mini.pw.edu.pl}
}

\corres{Piotr Graczyk}

\firstnote{submitted to MaxEnt 2022:  International Workshop on Bayesian Inference and Maximum Entropy Methods in Science and Engineering, IHP, Paris, July 18-22, 2022.} 

\abstract{We consider multivariate centered Gaussian models for the random vector $(Z^1,\ldots, Z^p)$, whose conditional structure is described by a homogeneous graph and which is  invariant under the action of a permutation subgroup. The following paper concerns with model selection within colored graphical Gaussian models, when the underlying conditional dependency graph is known.  We derive an analytic expression of the normalizing constant of the Diaconis–Ylvisaker conjugate  prior for the precision parameter and perform Bayesian  model selection in the class of graphical Gaussian models invariant by the action of a permutation subgroup. We illustrate our results with a toy example of dimension $5$.
}

\keyword{\textls[-15]{Graphical models; Colored graphical models; Invariance; Permutation symmetry; Diaconis-Ylvisaker conjugate prior }}

\begin{document}


\section{Introduction}
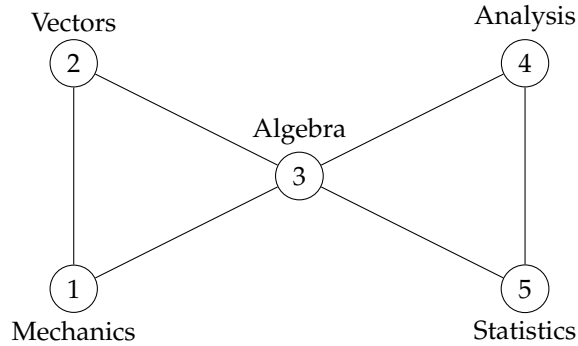
\begin{figure}[htb!]
\begin{center}
	\begin{tikzpicture}
		\node (A) [draw, circle] at (0,-3) {1};
		\draw (0,-3.3) node[below]{Mechanics};
		\node (B) [draw, circle] at (0, 0) {2};
		\draw (0, 0.3) node[above]{Vectors};
		\node (C) [draw, circle] at (3,-1.5) {3};
		\draw (3, -1.2) node[above]{Algebra};
		\node (D) [draw, circle] at (6, 0) {4};
		\draw (6, 0.3) node[above]{Analysis};
		\node (E) [draw, circle] at (6, -3) {5};
		\draw (6, -3.3) node[below]{Statistics};
		\draw (A) -- (B);
		\draw (A) -- (C);
		\draw (B) -- (C);
		\draw (C) -- (D);
		\draw (C) -- (E);
		\draw (D) -- (E);
	\end{tikzpicture}
\end{center}
	\caption{Conditional independence structure of examination marks.}
	\label{fig01}
\end{figure}

In the Graphical Gaussian model,
conditional independencies among components of 
a random vector $Z = (Z^1, Z^2, \ldots, Z^p)$
obeying the multivariate centered Gaussian law $\mathrm{N}(0, \Sigma)$ 
with an unknown covariance matrix 
$\Sigma \in \mathrm{Sym}^+(p, \real)$ 
are assigned by a simple undirected graph $\mathcal{G} = (\mathcal{V}, \mathcal{E})$,
where the set $\mathcal{V}$ of vertices is enumerated as $\mathcal{V} = \{1,2, \ldots, p\}$.
Namely, if the vertices $i$ and $j$ are disconnected in the graph $\mathcal{G}$,
then $Z^i$ and $Z^j$ are conditionally independent given other components $Z^k$, $k \ne i,j$.
This property is equivalent to that the $(i,j)$-component of the precision matrix $K := \Sigma^{-1}$ equals $0$.
Following \citet{HL08}, we impose the invariance on such a statistical model 
under the natural action of a permutation subgroup $\Gamma \subset \mathfrak{S}_p$
preserving the conditional independence structure of the model,
which means that
$\Gamma$ is a subgroup of the automorphism group 
$\mathrm{Aut}(\mathcal{G}) := \set{\sigma \in \mathfrak{S}_p}{\sigma(i) \sim \sigma(j)
\mbox{ if and only if } i \sim j} $
of the graph $\mathcal{G}$,
where $i \sim j$ means that there exists an edge between the vertices $i$ and $j$.
Such models are called RCOP graphical models.
It is proved that, when the graph $\mathcal{G}$ is homogeneous (i.e. decomposable and $A_4$-free, see \cite{LH07}),
the parameter set $\mathcal{P}_{\mathcal{G}}^{\Gamma}$ 
of precision matrices $K$ of our invariant model forms a homogeneous cone.
Therefore
we can apply our previous results \cite{GIK} about the Wishart laws on homogeneous cones
to this situation. 
In particular,
we obtain an exact analytic formula for the normalizing constant of 
the Diaconis-Ylvisaker conjugate prior of the precision matrix.
In order to demonstrate our results,
we work on the data set of 
the examination marks of $88$ students in $5$ different mathematical subjects
reported in \citet{MKB79}, following \citet{HL08}.
As is discussed in \citet{Wh90} and \citet{Ed00},
the data fit to the graphical Gaussian model from Fig. \ref{fig01}.

We carry out Bayesian model selection in the spirit of \citet{GrIsKoMa2022} of the group $\Gamma$ 
having the highest posterior probability among the ten possible groups 
preserving the graph above.   We note that in \cite{GrIsKoMa2022} only complete graph $\mathcal{G}$ was allowed. Thanks to new formulas for the normalizing constant of Diaconis-Ylvisaker conjugate prior, we are able to generalize results of \cite{GrIsKoMa2022} to homogeneous graphs.

The authors are grateful to anonymous referees for their careful readings and valuable comments.

\section{Main results}
Let us describe our results in more detail.
Let $\mathcal{Z}_{\mathcal{G}}^\Gamma$ be the linear space consisting of symmetric matrices $K \in \mathrm{Sym}(p, \real)$
such that $K_{\sigma(i) \sigma(j)} = K_{ij}$ for all $\sigma \in \Gamma$ and $i, j \in \mathcal{V}$, 
and $K_{ij} = 0$ if $i \ne j$ and $i \not\sim j$.
Then the cone $\mathcal{P}_{\mathcal{G}}^\Gamma$ equals
$\mathcal{Z}_{\mathcal{G}}^\Gamma \cap \mathrm{Sym}^+(p, \real)$,
so that our statistical model is the family of $N(0, \Sigma)$ with $\Sigma^{-1} \in \mathcal{P}_{\mathcal{G}}^\Gamma$.
The Diaconis-Ylvisaker conjugate prior for $K = \Sigma^{-1} \in \mathcal{P}_{\mathcal{G}}^\Gamma$ is given by
\begin{align}\label{eq:DY}
f(K;\delta, D) := \frac{1}{I_\mathcal{G}^\Gamma(\delta, D)} e^{- \mathrm{tr}\, KD /2}(\det K)^{(\delta-2)/2}  1_{\mathcal{P}_{\mathcal{G}}^\Gamma}(K)
\end{align}
for hyperparameters $\delta >2$ and $D \in \mathrm{Sym}^+(p, \real)$,
where
$$I_\mathcal{G}^\Gamma(\delta, D) := \int_{\mathcal{P}_{\mathcal{G}}^\Gamma} e^{- \mathrm{tr}\, KD /2}(\det K)^{(\delta-2)/2} \,dK$$
is the normalizing constant.
As is already stated, the cone $\mathcal{P}_{\mathcal{G}}^\Gamma$ is homogeneous, 
which means that there exists a linear group $H \subset GL(\mathcal{Z}_{\mathcal{G}}^\Gamma)$ acting on $\mathcal{P}_{\mathcal{G}}^\Gamma$ transitively.
Then, making use of our integral formula over $\mathcal{P}_{\mathcal{G}}^\Gamma$, see \eqref{eqn:factorization},
we can compute the normalizing constant $I_\mathcal{G}^\Gamma(\delta, D)$.

In order to state the integral formula, we introduce some functions.
Let $\Zsp$ be a linear subspace of $\mathrm{Sym}(p, \real)$ 
such that 
$\Pcone_{\Zsp} := \Zsp \cap \mathrm{Sym}^+(p, \real)$ is non-empty.
Let $\pi_\Zsp\colon \mathrm{Sym}(p, \real) \to \Zsp$ denote the orthogonal projection 
with respect to the trace inner product 
$\cpling{x}{y} := \mathrm{tr}\,xy$, $x, y \in \mathrm{Sym}(p, \real)$, that is,
\[
\cpling{x}{y}=\cpling{x}{\pi_\Zsp(y)},\qquad x\in\Zsp,\,\,y\in \mathrm{Sym}(p, \real).
\]
Let $\Pcone^*_\Zsp$ be the dual cone of $\Pcone_\Zsp$, that is,
$
\Pcone^*_\Zsp 
\colon= \set{y \in \Zsp}{ \cpling{x}{y} > 0 \mbox{ for all }x \in \lbar{\Pcone_\Zsp} \setminus \{0\} }.
$
It is easy to see that, if $D \in \mathrm{Sym}^+(p, \real)$, then $\pi_\Zsp(D) \in \Pcone_\Zsp^*$.
One can show
 that (see the proof of Proposition V.8 in \cite{Le14}),
for each $y \in \Pcone_\Zsp^*$, 
there exists a unique $\psi_\Zsp(y) \in \Pcone_\Zsp$ such that
the function $\Pcone_\Zsp\ni x\mapsto e^{-\cpling{x}{y}}\det x$ attains its maximum value at
$x = \psi_\Zsp(y)$, 
and that 
the map $\psi_\Zsp \colon \Pcone_\Zsp^* \to \Pcone_\Zsp$
equals the inverse map of
$\Pcone_\Zsp \owns x \mapsto \pi_\Zsp (x^{-1}) \in \Pcone_\Zsp^*$. 
For $y \in \Pcone^*_\Zsp$, 
define $\delta_\Zsp(y) := (\det \psi_\Zsp(y))^{-1}$ and
let $S_\Zsp(y) \colon \Zsp \to \Zsp$ be a linear operator defined
in such a way that
$$
\cpling{S_\Zsp(y) u}{v} = -\Bigl( \frac{\partial^2}{\partial s \partial t} \Bigr)
\log \delta_\Zsp(y + su + tv)\Big|_{s=t=0}
\qquad u,v \in \Zsp.
$$
Namely, $S_\Zsp(y)$ is the Hessian operator of a strictly convex function
$- \log \delta_\Zsp(y)$. 
Put $\varphi_\Zsp(y) := (\det S_\Zsp(y))^{1/2}$ for $y \in \Pcone^*_\Zsp$.
Finally, 
define 
$\gamma_\Zsp(\alpha) := \int_{\Pcone_\Zsp} e^{-\mathrm{tr}\,x} (\det x)^{\alpha}\, dx$
for $\alpha \ge 0$,
where $dx$ denotes the Lebesgue measure on $\Zsp$ normalized by the trace inner product.
Namely, $dx = \prod_{i=1}^{\dim \Zsp }dx_i$, where $(x_1, \ldots, x_{\dim \Zsp})$ is the standard coordinate system associated to an orthonormal basis of $\Zsp$ with respect to the trace inner product.

\begin{Theorem} \label{thm:mainthm}
	If $\Zsp = \Zsp_{\mathcal{G}}^\Gamma$, then one has
	\begin{equation}  \label{eqn:factorization}
		\int_{\Pcone_\Zsp} e^{- \cpling{x}{y}} (\det x)^\alpha \,dx 
		= \gamma_\Zsp(\alpha)\, \varphi_\Zsp(y) \,\delta_\Zsp(y)^{-\alpha}
		\qquad y \in \Pcone_\Zsp^*,\,\,\alpha \ge 0.
	\end{equation}
\end{Theorem}

We shall show Theorem~\ref{thm:mainthm} by using the homogeneity of
$\Pcone_{\mathcal{G}}^\Gamma =\Pcone_{\Zsp_\mathcal{G}^\Gamma}$
in our case,
whereas we notice that the formula (\ref{eqn:factorization}) is valid also 
for some non-homogeneous cases, 
e.g. the cone $\Pcone_\Zsp$ arising from uncolored decomposable graphical models \cite{Ro00}.

Let $\gamma_{\mathcal{G}}^\Gamma, \,\delta_{\mathcal{G}}^\Gamma$
and $\varphi_{\mathcal{G}}^\Gamma$
denote the functions $\gamma_\Zsp,\,\, \delta_\Zsp$ and $\varphi_\Zsp$ respectively
with $\Zsp = \Zsp_{\mathcal{G}}^\Gamma$.
Then we have
$$ I_\mathcal{G}^\Gamma(\delta, D) =
\gamma_\mathcal{G}^\Gamma((\delta-2)/2)\,
\varphi_\mathcal{G}^\Gamma(\pi_\Zsp(D)/2)\,
\delta_\mathcal{G}^\Gamma(\pi_\Zsp(D)/2)^{-(\delta -2)/2}. $$
In our Bayesian model selection setting \cite{GrIsKoMa2022},
for a fixed graph $\mathcal{G}$, we suppose that a group $\Gamma$
is distributed uniformly over all the possible subgroups of the automorphism group of
the graph $\mathcal{G}$.
Given samples $Z_1,\, Z_2, \dots, Z_n$, that is, the i.i.d.~random vectors
obeying $N(0, \Sigma)$ with $\Sigma^{-1} \in \Pcone_{\mathcal{G}}^\Gamma$,
we see that the posterior probability $\mathbb{P}(\Gamma | Z_1, \ldots, Z_n)$ 
is proportional to
$ I_{\mathcal{G}}^\Gamma(\delta + n, D + \sum_{i=1}^n Z_i \transp{Z_i}) 
/ I_{\mathcal{G}}^\Gamma(\delta, D)$, see \cite[Eq. (30)]{GrIsKoMa2022}.

\section{Matrix realization of homogeneous cones}
It is known that any homogeneous cone is linearly isomorphic to
some $\Pcone_{\Zsp}$, 
where $\Zsp \subset \mathrm{Sym}(p, \real)$ is a linear subspace
consisting of real symmetric matrices 
admitting certain specific block decomposition described below
\cite{Is14, Is15}. 
Let $n_1, n_2, \ldots, n_r$ be positive integers such that
$p= n_1 + n_2 + \ldots + n_r$. 
Let $V = \{ V_{lk} \}_{1 \le k < l \le r}$
be a family of linear spaces $V_{lk}\subset \mathrm{Mat}(n_l, n_k; \real)$
satisfying the following axioms:\\
(V1) $A \in V_{lk} \Rightarrow A \transp{A} \in \real I_{n_l}
\,\,\,(1 \le k < l \le r)$,\\
(V2) $A \in V_{lj},\, B \in V_{kj} \Rightarrow A \transp{B} \in V_{lk}
\,\,\,(1 \le j < k < l \le r)$,\\
(V3) $A \in V_{lk},\, B \in V_{kj} \Rightarrow AB \in V_{lj}
\,\,\,(1 \le j < k < l \le r)$.\\
Let $\Zsp_V$ be the linear space consisting of $x \in \mathrm{Sym}(p, \real)$
of the following form:
\begin{equation} \label{eqn:x}
	x = \begin{pmatrix} X_{11} & \transp{X_{21}} & \cdots &\transp{X_{r1}} \\
		X_{21} & X_{22} & & \transp{X_{r2}} \\
		\vdots & & \ddots &  \vdots \\
		X_{r1} & X_{r2} & \cdots & X_{rr} \end{pmatrix}
	\quad
	\begin{aligned}
		X_{kk} &= x_{kk} I_{n_k},\,x_{kk} \in \real,\,\,k=1, \ldots, r,\\
		X_{lk} &\in V_{lk},\,\,1 \le k < l \le r.
	\end{aligned}
\end{equation} 
Let $H_V$ be the set of $p \times p$ lower triangular matrices $T$ of the following form:
\begin{equation} \label{eqn:T}
	T = \begin{pmatrix} T_{11} & & & \\ T_{21} & T_{22} & & \\ \vdots & & \ddots & \\
		T_{r1} & T_{r2} & \cdots & T_{rr} \end{pmatrix}
	\quad
	\begin{aligned}
		T_{kk} &= t_{kk} I_{n_k},\,t_{kk}>0,\,\,k=1, \ldots, r,\\
		T_{lk} &\in V_{lk},\,\,1 \le k < l \le r.
	\end{aligned}
\end{equation}
Then $H_V$ forms a Lie group by (V3), 
and acts on $\Zsp_V$ linearly by
$\rho(T)x := T x \transp{T}$, $T \in H_V$, $x \in \Zsp_V$.
Moreover, the group $H_V$ acts on the cone $\Pcone_V := \Pcone_{\Zsp_V}$
simply transitively by $\rho$.
We write $\rho^*(T)$, $T \in H_V$, for 
the adjoint of  the linear operator $\rho(T)$ on $\Zsp_V$ with respect to the trace inner product,
which means that
$\cpling{\rho(T)x}{y} = \cpling{x}{\rho^*(T)y}$ for $x, y \in \Zsp_V$.
Then we see that
$\rho^*(T)y = \pi_{\Zsp}(\transp{T} y T )$ for $y \in \Zsp_V$.
Moreover, 
if $y \in \Pcone^*_V := \Pcone^*_{\Zsp_V}$,
we have
$\psi_{\Zsp_V}(\rho^*(T)y) = \rho(T^{-1})\psi_{\Zsp_V}(y)$. 
Since $\psi_\Zsp(y)  \in \Pcone_V$ for $y \in \Pcone^*_V$,
we can take a unique $T_y \in H_V$ for which $\psi_\Zsp(y) = \rho(T_y^{-1}) I_p$. 
Then we have $y = \rho^*(T_y) I_p$
because
$\rho(T_y^{-1}) I_p = \rho(T_y^{-1})\psi_\Zsp(I_p) =  \psi_\Zsp( \rho^*(T_y) I_p)$. 

\begin{Lemma} \label{lemma:rel_inv}
	For $y \in \Pcone^*_{\Zsp_V}$,
	one has $\delta_{\Zsp_V}(y) = (\det T_y)^2$ 
	and $\varphi_{\Zsp_V}(y) = (\det \rho(T_y))^{-1}$.
\end{Lemma}
Using Lemma~\ref{lemma:rel_inv} and
a general theory about relatively invariant functions on a homogeneous cone
(see e.g. \cite[Section IV]{Is14}), 
we can compute explicitly $\delta_{\Zsp_V}$ and $\varphi_{\Zsp_V}$.

Put $q_k := \sum_{l>k} \dim V_{lk}$ for $k=1, \ldots, r$ and $N := \dim \Zsp_V$. 
\begin{Theorem} \label{thm:Gindikin_integral}
	{\rm (i)} One has
	$$\gamma_{\Zsp_V}(\alpha) 
	= (2\pi)^{(N-r)/2} 
	\prod_{k=1}^r
	\Bigl( n_k^{-n_k \alpha-(q_k +1)/2 } \Gamma(n_k \alpha + (q_k/2) + 1) \Bigr).$$
	{\rm (ii)} The equality {\rm (\ref{eqn:factorization})} holds for $\Zsp = \Zsp_V$. 
\end{Theorem}
We give a sketch of the proof of Theorem \ref{thm:Gindikin_integral}.
We denote by $(A|B)$ the trace inner product
$\mathrm{tr}\, A \transp{B}$ of $A,\, B \in V_{lk}$.
Then for an element $x \in \Zsp_V$ in (\ref{eqn:x}), we have
$$
\cpling{x}{x} = \sum_{k=1}^r (n_k x_{kk}^2) + 2 \sum_{1 \le k < l \le r} (X_{lk} | X_{lk}),
$$
so that
$
dx = \prod_{k=1}^r (n_k^{1/2} dx_{kk}) \prod_{1 \le k < l \le r} (2^{\dim V_{lk}/2} dX_{lk}),
$
where $dX_{lk}$ stands for the Lebesgue measure on $V_{lk}$
normalized by $(\cdot | \cdot)$.
By the change of variable $x = \rho(T)I_p$ with $T \in H_V$ in (\ref{eqn:T}),
we get
$$
dx = \prod_{k=1}^r (n_k^{1/2} 2 t_{kk}^{1 + q_k} dt_{kk}) 
\prod_{1 \le k < l \le r} 2^{\dim V_{lk}/2} dT_{lk},
$$
so that $\gamma_{\Zsp_V}(\alpha)$ equals
\begin{align*}
	\prod_{k=1}^r \Bigl(
	\int_0^\infty e^{- n_k t_{kk}^2} n_k^{1/2} 2 t_{kk}^{2n_k \alpha + 1 + q_k}\,dt_{kk} 
	\Bigr)
	\prod_{1 \le k < l \le r} \Bigl(2^{\dim V_{lk}/2} \int_{V_{lk}} e^{- (T_{lk} | T_{lk})}\,dT_{lk}\Bigr). 
\end{align*}
As for (ii), we observe that
$ \cpling{x}{y} = \cpling{x}{\rho^*(T_y)I_p} = \cpling{\rho(T_y)x}{I_p} =\mathrm{tr}\, \rho(T_y)x.$
By the change of variable $x' = \rho(T_y) x$,
we have
\begin{align*}
	\int_{\Pcone_V} e^{- \cpling{x}{y}} (\det x)^\alpha \,dx
	&= \int_{\Pcone_V} e^{- \mathrm{tr}\, x'} \{(\det x')(\det T_y)^{-2} \}^\alpha\,
	(\det \rho(T_y))^{-1} dx'\\
	&= (\det T_y)^{- 2\alpha} (\det \rho(T_y))^{-1} \gamma_{\Zsp_V}(\alpha),
\end{align*}  
so that (\ref{eqn:factorization}) follows from Lemma \ref{lemma:rel_inv}.

\begin{Theorem} \label{thm:conjugation}
	For a homogeneous graph $\mathcal{G} = (\mathcal{V},\, \mathcal{E})$
	and a subgroup $\Gamma$ of the automorphism group of $\mathcal{G}$,
	there exists an orthogonal matrix $U \in O(p)$ such that
	$\transp{U} \Zsp_{\mathcal{G}}^\Gamma U = \Zsp_V$
	with some $V = \{ V_{lk} \}_{1 \le k < l \le r}$.
\end{Theorem}
The proof is omitted, it is based on a representation theory similarly as in \cite{GrIsKoMa2022} and uses a proper ordering of vertices.
Theorem \ref{thm:conjugation} together with Theorem~\ref{thm:Gindikin_integral} yields
Theorem~\ref{thm:mainthm}.  Theorem \ref{thm:conjugation} is very important from the practical point of view. The knowledge of the orthogonal matrix $U$ allows to identify all parameters of the space $\Zsp_V$ (see \eqref{eqn:x}). However, the problem of finding a suitable $U$ matrix is in general very complicated. In the next section we will consider the butterfly model  from Fig. \ref{fig01} and present the exact forms of the $U$ matrices for all subgroups of $\mathrm{Aut}(\mathcal{G})$.

\section{Toy example}
In what follows, let $\mathcal{G}$ be the five-vertex graph from Fig. \ref{fig01}.
We use the cyclic representation of permutations on $\mathcal{V}=\{1,\ldots,5\}$.
Then the group $\mathrm{Aut}(\mathcal{G})$ is generated by
$\sigma_1 := \begin{pmatrix} 1 & 2 \end{pmatrix},\,
\sigma_2 := \begin{pmatrix} 4 & 5 \end{pmatrix}$,\,
and $\sigma_3 := \begin{pmatrix} 1 & 4 \end{pmatrix} \begin{pmatrix} 2 & 5 \end{pmatrix}$.
Put
$\tau := \sigma_2 \circ \sigma_3 = \begin{pmatrix} 1 & 5 & 2 & 4 \end{pmatrix}$.
Then 
$\sigma_2 = \sigma_1 \circ \tau^2$ and $\sigma_3 = \sigma_1 \circ \tau^3$,
so that $\mathrm{Aut}(\mathcal{G})$ 
is generated by $\sigma_1$ and $\tau$.
Moreover, since the orders of $\sigma_1$ and $\tau$ are $2$ and $4$ respectively
with $\sigma_1 \circ \tau \circ \sigma_1^{-1} =\tau^3$,
the group $\mathrm{Aut}(\mathcal{G})$ 
equals the dihedral group 
$\langle \sigma_1 \rangle \ltimes \langle \tau \rangle$
 of order $8$.   
Then all the subgroups of $\mathrm{Aut}(\mathcal{G})$ are listed as:
\begin{gather*}
	\Gamma_1 := \{ e\}, \quad 
	\Gamma_2 := \langle \sigma_1 \rangle, \quad
	\Gamma_3 := \langle \sigma_1 \circ \tau \rangle, \quad
	\Gamma_4 := \langle \sigma_1 \circ \tau^2\rangle, \quad
	\Gamma_5 := \langle \sigma_1 \circ \tau^3\rangle,\\
	\Gamma_6 := \langle \tau^2 \rangle, \quad
	\Gamma_7 := \langle \tau \rangle, \quad
	\Gamma_8 := \langle \sigma_1,\, \tau^2 \rangle, \quad
	\Gamma_9 := \langle \sigma_1\circ \tau,\, \tau^2  \rangle, \quad
	\Gamma_{10} := \mathrm{Aut}(\mathcal{G}).
\end{gather*}
(i) When $\Gamma = \Gamma_1$, then $x \in \Zsp_{\mathcal{G}}^\Gamma$ is of the form
$$
x =\begin{pmatrix} x_{11} & x_{21} & x_{31} & 0 & 0\\
	x_{21} & x_{22} & x_{32} & 0 & 0 \\
	x_{31} & x_{32} & x_{33} & x_{43} & x_{53} \\
	0 & 0 & x_{43} & x_{44} & x_{54} \\
	0 & 0 & x_{53} & x_{54} & x_{55} \end{pmatrix}.
$$
We give an orthogonal matrix $U$ so that $\transp{U} x U$ becomes of the form of matrix realization in the previous section (see Theorem \ref{thm:conjugation}) by
$$
U := \begin{pmatrix} 1 & 0 & 0 & 0 & 0\\ 0 & 1 & 0 & 0 & 0 \\ 0 & 0 & 0 & 0 & 1 \\
	0 & 0 & 1 & 0 & 0 \\ 0 & 0 & 0 & 1 & 0 \end{pmatrix},\quad
\transp{U} x U 
= \begin{pmatrix} x_{11} & x_{21} & 0& 0 & x_{31} \\
	x_{21} & x_{22} & 0 & 0 & x_{32}\\
	0 & 0 & x_{44} & x_{54} & x_{43}\\
	0 & 0 & x_{54} & x_{55} & x_{53}\\
	x_{31} & x_{32} & x_{43} & x_{53} & x_{33} \end{pmatrix}.
$$
In this case $N=11,\, r=5,\, n_1 = \ldots = n_5 = 1$, and $V_{lk}\,\,\,(1 \le k < l \le 5)$ is $\mathbb{R}$ or $\{0\}$.
Since $\gamma_{\transp{U}\Zsp U}(\alpha) = \gamma_{\Zsp}(\alpha)$ in general,
we see from Theorem~\ref{thm:Gindikin_integral} (i) that  
$$
\gamma_{\mathcal{G}}^\Gamma(\alpha) 
= (2\pi)^3 \Gamma(\alpha +1)\Gamma(\alpha + \frac{3}{2})^2 \Gamma(\alpha + 2)^2.$$
Moreover, 
the functions 
$\delta_{\mathcal{G}}^\Gamma(x)$ and
$\varphi_{\mathcal{G}}^\Gamma(x)$ are expressed respectively as 
$$
x_{33}^{-1} 
\begin{vmatrix} x_{11} & x_{21} & x_{31} \\ 
	x_{21} & x_{22} & x_{32} \\
	x_{31} & x_{32} & x_{33} \end{vmatrix}
\begin{vmatrix} x_{33} & x_{43} & x_{53} \\
	x_{43} & x_{44} & x_{54} \\
	x_{53} & x_{54} & x_{55} \end{vmatrix}, \quad
x_{33} 
\begin{vmatrix} x_{11} & x_{21} & x_{31} \\ 
	x_{21} & x_{22} & x_{32} \\
	x_{31} & x_{32} & x_{33} \end{vmatrix}^{-2}
\begin{vmatrix} x_{33} & x_{43} & x_{53} \\
	x_{43} & x_{44} & x_{54} \\
	x_{53} & x_{54} & x_{55} \end{vmatrix}^{-2}.
$$
For $i=2, \ldots, 10$, 
clearly $\Pcone_{\mathcal{G}}^{\Gamma_i}$
is a subset of 
$\Pcone_{\mathcal{G}}^{\Gamma} = \Pcone_{\mathcal{G}}^{\Gamma_1}$
and one can show that
$\delta_{\mathcal{G}}^{\Gamma_i}$ is equal to the restriction of 
$\delta_{\mathcal{G}}^{\Gamma}$ above to $\Pcone_{\mathcal{G}}^{\Gamma_i}$.\\

\noindent
(ii) When $\Gamma = \Gamma_2$, 
we describe respectively 
$x \in \Zsp_{\mathcal{G}}^{\Gamma}$,
an orthogonal matrix $U$ and $\transp{U} x U$ as
$$
\begin{pmatrix} a & b & c & 0& 0 \\
	b & a & c & 0 & 0 \\
	c & c & d & e & f \\
	0 & 0 & e & g & h \\
	0 & 0 & f & h & i \end{pmatrix}, 
\quad
\begin{pmatrix} 1/\sqrt{2} & 1/\sqrt{2} & 0 & 0 & 0 \\
	-1/\sqrt{2} & 1/\sqrt{2} & 0 & 0 & 0 \\
	0 & 0 & 0 & 0 & 1 \\
	0 & 0 & 1 & 0 & 0\\
	0 & 0 & 0 & 1 & 0 \end{pmatrix}, 
\quad
\begin{pmatrix} a-b & 0 & 0 & 0 & 0\\
	0 & a+b & 0 & 0 & \sqrt{2} c \\
	0 & 0 & g & h & e \\
	0 & 0 & h & i & f \\
	0 & \sqrt{2} c & e & f & d \end{pmatrix}.
$$
We have
$\gamma_{\mathcal{G}}^\Gamma(\alpha)
= (2\pi)^2 \Gamma(\alpha+1)^2 \Gamma(\alpha + \frac{3}{2})^2 \Gamma(\alpha+2)$,
and
$$\varphi_{\mathcal{G}}^\Gamma(x)
=d (a-b)^{-1} 
\begin{vmatrix} a+b & \sqrt{2} c \\ \sqrt{2} c & d\end{vmatrix}^{-3/2}
\begin{vmatrix} d &e & f \\ e & g  & h \\ f & h & i \end{vmatrix}^{-2}. $$

\noindent 
(iii) When $\Gamma = \Gamma_3$, 
we describe respectively $x \in \Zsp_{\mathcal{G}}^\Gamma$,
an orthogonal matrix $U$ and $\transp{U} x U$ as
$$
\begin{pmatrix} a & d & e & 0 & 0 \\
	d & b & f & 0 & 0 \\
	e & f & c & f & e \\
	0 & 0 & f & b & d \\
	0 & 0 & e & d & a \end{pmatrix}, 
\quad
\begin{pmatrix} 1 & 0 & 0 & 0 & 0 \\
	0 & 0 & 1 & 0 & 0\\
	0 & 0 & 0 & 0 & 1\\
	0 & 0 & 0 & 1 & 0\\
	0 & 1 & 0 & 0 & 0\end{pmatrix},
\quad
\begin{pmatrix} a & 0 & d & 0 & e\\
	0 & a & 0 & d & e\\
	d & 0 & b & 0 & f\\
	0 & d & 0& b & f\\
	e & e & f & f & c\end{pmatrix}.
$$
Then
$\gamma_{\mathcal{G}}^\Gamma(\alpha)
= (2\pi)^{3/2} 2^{-4\alpha-5/2} \Gamma(2\alpha+2) \Gamma(2\alpha + \frac{3}{2}) \Gamma(\alpha+1)$
and
$\varphi_{\mathcal{G}}^\Gamma(x) =
\begin{vmatrix} a & d & e\\ d & b & f \\ e & f & c \end{vmatrix}^{-2}.$\\

\noindent 
(iv) When $\Gamma = \Gamma_4$, 
we describe respectively $x \in \Zsp_{\mathcal{G}}^\Gamma$,
an orthogonal matrix $U$ and $\transp{U} x U$ as
$$
\begin{pmatrix} a & d & e & 0 & 0 \\
	d & b & f & 0 & 0 \\
	e & f & c & g & g \\
	0 & 0 & g & h & i \\
	0 & 0 & g & i & h \end{pmatrix}, 
\quad
\begin{pmatrix} 1 & 0 & 0 & 0 & 0 \\
	0 & 1 & 0 & 0 & 0\\
	0 & 0 & 0 & 0 & 1\\
	0 & 0 & 1/\sqrt{2} & 1/\sqrt{2} & 0\\
	0 & 0 & -1/\sqrt{2}& 1/\sqrt{2} & 0 \end{pmatrix},
\quad
\begin{pmatrix} a & d & 0 & 0 & e \\
	d & b & 0 & 0 & f \\
	0 & 0 & h-i & 0 & 0 \\
	0 & 0 &  0 & h+i& \sqrt{2}g \\
	e & f & 0 & \sqrt{2}g & c \end{pmatrix}.
$$
Then 
$\gamma_{\mathcal{G}}^\Gamma(\alpha)
= \gamma_{\mathcal{G}}^{\Gamma_2}(\alpha)
= (2\pi)^2 \Gamma(\alpha+1)^2 \Gamma(\alpha + \frac{3}{2})^2 \Gamma(\alpha+2)$,
while $\varphi_{\mathcal{G}}^\Gamma(x)$ is equal to 
$$c (h-i)^{-1} 
\begin{vmatrix} h+i & \sqrt{2} g \\ \sqrt{2} g & c\end{vmatrix}^{-3/2}
\begin{vmatrix} a &d & e \\ d & b  & f \\ e & f & c \end{vmatrix}^{-2}. $$

\noindent 
(v) When $\Gamma = \Gamma_5$, 
we describe respectively $x \in \Zsp_{\mathcal{G}}^\Gamma$,
an orthogonal matrix $U$ and $\transp{U} x U$ as
$$
\begin{pmatrix} a & d & e & 0 & 0 \\
	d & b & f & 0 & 0 \\
	e & f & c & e & f \\
	0 & 0 & e & a & d\\
	0 & 0 & f & d & b \end{pmatrix},
\quad
\begin{pmatrix} 1 & 0 & 0 & 0 & 0\\
	0 & 0 & 1 & 0 & 0\\
	0 & 0 & 0 & 0 & 1\\
	0 & 1 & 0 & 0 & 0\\
	0 & 0 & 0 & 1 & 0\end{pmatrix},
\quad
\begin{pmatrix} a & 0 & d & 0 & e\\
	0 & a & 0 & d & e\\
	d & 0 & b & 0 & f\\
	0 & d & 0& b & f\\
	e & e & f & f & c\end{pmatrix}.
$$
Then
$\gamma_{\mathcal{G}}^\Gamma(\alpha) = \gamma_{\mathcal{G}}^{\Gamma_3}(\alpha) 
= (2\pi)^{3/2} 2^{-4\alpha-5/2} \Gamma(2\alpha+2) \Gamma(2\alpha + \frac{3}{2}) \Gamma(\alpha+1)$
and
$\varphi_{\mathcal{G}}^\Gamma(x) =
\begin{vmatrix} a & d & e\\ d & b & f \\ e & f & c \end{vmatrix}^{-2}.$\\
(vi) We have $\Zsp_{\mathcal{G}}^{\Gamma_6}=\Zsp_{\mathcal{G}}^{\Gamma_8}$. When $\Gamma = \Gamma_6$ or $\Gamma_8$, 
we describe respectively $x \in \Zsp_{\mathcal{G}}^\Gamma$,
an orthogonal matrix $U$ and $\transp{U} x U$ as
\begin{gather*}
	\begin{pmatrix} a & b & c & 0 & 0 \\
		b & a & c & 0 & 0 \\
		c & c & d & e & e \\
		0 & 0 & e & f & g\\
		0 & 0 & e & g & f \end{pmatrix},
	\quad
	\begin{pmatrix} 1/\sqrt{2} & 1/\sqrt{2} & 0  & 0 & 0 \\
		-1/\sqrt{2} & 1/\sqrt{2} & 0 & 0 & 0\\
		0 & 0 & 0 & 0 & 1\\
		0 & 0 & 1/\sqrt{2} & 1/\sqrt{2} & 0 \\
		0 & 0 &-1/\sqrt{2} & 1/\sqrt{2} & 0\end{pmatrix},\\
	\begin{pmatrix} a-b & 0 & 0 & 0 & 0 \\
		0 & a+b & 0 & 0 & \sqrt{2} c\\
		0 & 0 & f-g & 0 & 0\\
		0 & 0 & 0 & f+g & \sqrt{2} e\\
		0 & \sqrt{2} c& 0 & \sqrt{2} e & d \end{pmatrix}.
\end{gather*}
Then
$\gamma_{\mathcal{G}}^\Gamma(\alpha) =
2\pi \Gamma(\alpha+1)^3 \Gamma(\alpha + \frac{3}{2})^2$,
and
$$ \varphi_{\mathcal{G}}^\Gamma(x)
= d (a-b)^{-1}(f-g)^{-1}
\begin{vmatrix} a+b & \sqrt{2} c\\ \sqrt{2} c & d \end{vmatrix}^{-3/2}
\begin{vmatrix} f+g & \sqrt{2} e\\ \sqrt{2}e & d\end{vmatrix}^{-3/2}.$$

\noindent
(vii) For $\Gamma = \Gamma_7,\, \Gamma_9$ or $\Gamma_{10}$  the linear space $\Zsp_{\mathcal{G}}^\Gamma$ is the same. We describe respectively $x \in \Zsp_{\mathcal{G}}^\Gamma$,
an orthogonal matrix $U$ and $\transp{U} x U$ as
\begin{gather*}
	\begin{pmatrix} a & c & d & 0 & 0\\
		c & a & d & 0 & 0\\
		d & d & b & d & d\\
		0 & 0 & d & a & c\\
		0 & 0 & d & c & a\end{pmatrix},
	\quad
	\begin{pmatrix} 1/\sqrt{2} & 0 & 1/\sqrt{2} & 0 & 0\\
		-1/\sqrt{2} & 0 & 1/\sqrt{2} & 0 & 0\\
		0 & 0 & 0 & 0 &1\\
		0 & -1/\sqrt{2} & 0 & 1/\sqrt{2} & 0\\
		0 & 1/\sqrt{2} & 0 & 1/\sqrt{2} & 0\end{pmatrix},\\
	\begin{pmatrix} a-c & 0 & 0 & 0 & 0\\
		0 & a-c & 0 & 0 & 0\\
		0 & 0 & a+c & 0 & \sqrt{2}d\\
		0 & 0 & 0 & a+c & \sqrt{2}d\\
		0 & 0 & \sqrt{2}d & \sqrt{2}d & b \end{pmatrix}.
\end{gather*}
Then
$\gamma_{\mathcal{G}}^\Gamma(\alpha) =
(2\pi)^{1/2} {2^{-4\alpha - 3/2}} \Gamma(2\alpha+1)\Gamma(2\alpha+\frac{3}{2})\Gamma(\alpha+1)$
and 
$ \varphi_{\mathcal{G}}^\Gamma(x)
= (a-c)^{-1}\begin{vmatrix} a+c &  \sqrt{2} d\\ \sqrt{2}d & b \end{vmatrix}^{-3/2}.$

\section{Numerical example}
We carry out Bayesian model selection in the spirit of \citet{GrIsKoMa2022}.
For a fixed graph $\mathcal{G}$, we suppose that a group $\Gamma$
is distributed uniformly over all possible subgroups of $\mathrm{Aut}(\mathcal{G})$. Given sample $Z_1,\ldots, Z_n$ from $\mathrm{N}(0, \Sigma)$ with $\Sigma^{-1} \in \Pcone_{\mathcal{G}}^\Gamma$, where $K=\Sigma^{-1}$ follows the Diaconis-Ylvisaker conjugate prior \eqref{eq:DY} with hyperparameters $(\delta,D)$, then the posterior probability $\mathbb{P}(\Gamma | Z_1, \ldots, Z_n)$ 
is proportional to 
$ I_{\mathcal{G}}^\Gamma(\delta + n, D + \sum_{i=1}^n Z_i \transp{Z_i}) 
/ I_{\mathcal{G}}^\Gamma(\delta, D)$.
In order to demonstrate our results,
we work on the data set of 
the examination marks of $n=88$ students in $p=5$ different mathematical subjects. As was reported in \cite{Wh90} and  \cite{Ed00}, the data demonstrate an excellent fit to the graphical Gaussian model displayed in Fig. \ref{fig01}.
Since the groups $\Gamma=\Gamma_6, \Gamma_8$ (similarly $\Gamma_7$, $\Gamma_9$ and $\Gamma_{10}$) impose the same symmetries on $\Zsp_{\mathcal{G}}^\Gamma$, we consider Bayesian model selection within $7$ different models: $\Gamma_i$ for $i=1, 2,\ldots, 7$.
The five mathematical subjects are enumerated as in the graph presented in Fig. \ref{fig01}. 
As our method applies only to centered normal sample, as usual, we center the marks and consider a correction of the degrees of freedom $n^*=88-1=87$. Then,
\[
\sum_{i=1}^n Z_i \transp{Z_i} = \begin{pmatrix}
	26\,601.82 & 11\,068.36 &  \ \ 8\,837.41 &  \ \ 9\,245.73 & 10\,214.23 \\
	11\,068.36 & 15\,037.27 &  \ \ 7\,408.68 &  \ \ 8\,236.55 &  \ \ |8\,614.05 \\
	  \ \ 8\,837.41 & \ \  7\,408.68 &  \ \ 9\,821.08 & \ \  9\,753.86 & 10\,602.74 \\
	\ \    9\,245.73 &  \ \ 8\,236.55 &\ \   9\,753.86 & 19\,173.09 & 13\,531.59 \\
	 10\,214.23 &  \ \ 8\,614.05 & 10\,602.74 & 13\,531.59 & 25\,904.72
\end{pmatrix}.
\]
We take usual hyperparameter $\delta=3$ and $D=d\cdot I_5$ for $d\in\{1,10^2,10^4\}$. Below we present a subgroup with the highest posterior probability $p$. 
\[
\begin{array}{c|c|c}
	 d=1 & d=100 & d=10\,000 \\
\hline
\Gamma_7 \,\,(p=1) & \Gamma_3 \,\,(p=0.8) & \Gamma_1 \,\,(p=0.75)
\end{array}
\]
Depending on the value of hyperparameter $D$, the model with highest posterior probability is
\begin{itemize}
	\item $\Gamma_7$, which corresponds to full symmetry as $\Zsp_{\mathcal{G}}^{\Gamma_7}=\Zsp_{\mathcal{G}}^{\mathrm{Aut}(\mathcal{G})}$,
	\item $\Gamma_3 = \left< \begin{pmatrix}
		1 & 5
	\end{pmatrix}\begin{pmatrix}
	2 & 4
\end{pmatrix}\right>$, which corresponds to invariance of the model to interchange (Mechanics, Vectors) $\leftrightarrow$ (Statistics, Analysis),
\item $\Gamma_1=\{e\}$, which corresponds to no additional symmetry.
\end{itemize}
The hyperparameter $\delta$ has much less impact on model selection.

We note that the same example was considered in \cite{HL08}, where the colored graphical models were introduced for the first time. The authors of \cite{HL08}, using a BIC criterion, point out that  model $\Gamma_3$  (see \cite[Fig. 8]{HL08}) represents an excellent fit. 

Fitted concentrations $\times 10^3$ for the examination marks assuming the model $\Gamma_3$ are presented below.\footnote{In \cite[Table 6]{HL08}  erroneous entries are presented in the same table.}
\[
\begin{array}{r|rrrrr}
	&  Mechanics &  Vectors &  Algebra &  Analysis & Statistics \\
	\hline
Mechanics	& 5.85 &  -2.23 & -3.72  & 0 &  0 \\
Vectors	& -2.23  & 10.15 & -5.88 & 0 & 0 \\
Algebra	& -3.72 & -5.88 & 26.95 & -5.88 & -3.72 \\
Analysis	& 0 & 0 & -5.88 & 10.15 & -2.23 \\
Statistics	& 0 & 0 & -3.72 & -2.23 & 5.85\\\hline
\end{array}
\]

\funding{Research of P. Graczyk was supported by JST PRESTO.
Research of H. Ishi was supported by JST PRESTO, KAKENHI 20K03657, and Osaka Central Advanced Mathematical Institute (MEXT Joint Usage/Research Center on Mathematics and Theoretical Physics JPMXP0619217849).
Research of B. Ko{\l}odziejek was funded by (POB Cybersecurity and Data Science) of Warsaw University of Technology within the Excellence Initiative: Research University (IDUB) programme.}

\reftitle{References}

\bibliography{biblio/CGM}

\begin{thebibliography}{-------}
\providecommand{\natexlab}[1]{#1}

\bibitem[H{\o}jsgaard and Lauritzen(2008)]{HL08}
H{\o}jsgaard, S.; Lauritzen, S.L.
\newblock Graphical {G}aussian models with edge and vertex symmetries.
\newblock {\em J. R. Stat. Soc. Ser. B Stat. Methodol.} {\bf 2008}, {\em
  70},~1005--1027.

\bibitem[Letac and Massam(2007)]{LH07}
Letac, G.; Massam, H.
\newblock Wishart distributions for decomposable graphs.
\newblock {\em Ann. Statist.} {\bf 2007}, {\em 35},~1278--1323.

\bibitem[Graczyk \em{et~al.}(2019)Graczyk, Ishi, and Ko{\l}odziejek]{GIK}
Graczyk, P.; Ishi, H.; Ko{\l}odziejek, B.
\newblock Wishart laws and variance function on homogeneous cones.
\newblock {\em Probab. Math. Statist.} {\bf 2019}, {\em 39},~337--360.

\bibitem[Mardia \em{et~al.}(1979)Mardia, Kent, and Bibby]{MKB79}
Mardia, K.V.; Kent, J.T.; Bibby, J.M.
\newblock {\em Multivariate analysis}; Probability and Mathematical Statistics:
  A Series of Monographs and Textbooks, Academic Press [Harcourt Brace
  Jovanovich, Publishers], London-New York-Toronto, Ont.,  1979; pp. xv+521.

\bibitem[Whittaker(1990)]{Wh90}
Whittaker, J.
\newblock {\em Graphical models in applied multivariate statistics}; Wiley
  Series in Probability and Mathematical Statistics: Probability and
  Mathematical Statistics, John Wiley \& Sons, Ltd., Chichester,  1990; pp.
  xiv+448.

\bibitem[Edwards(2000)]{Ed00}
Edwards, D.
\newblock {\em Introduction to graphical modelling}, second ed.; Springer Texts
  in Statistics, Springer-Verlag, New York,  2000; pp. xvi+333.

\bibitem[Graczyk \em{et~al.}(2022)Graczyk, Ishi, Ko{\l}odziejek, and
  Massam]{GrIsKoMa2022}
Graczyk, P.; Ishi, H.; Ko{\l}odziejek, B.; Massam, H.
\newblock Model selection in the space of Gaussian models invariant by
  symmetry.
\newblock {\em Ann. Statist.} {\bf 2022}, {\em 50},~1747--1774.

\bibitem[Letac(2014)]{Le14}
Letac, G.
\newblock Decomposable Graphs. In {\em Modern methods of multivariate
  statistics}; Hermann,  2014; Vol.~82, pp. 155--196.

\bibitem[Roverato(2000)]{Ro00}
Roverato, A.
\newblock Cholesky decomposition of a hyper inverse {W}ishart matrix.
\newblock {\em Biometrika} {\bf 2000}, {\em 87},~99--112.

\bibitem[Ishi(2014)]{Is14}
Ishi, H.
\newblock Homogeneous cones and their applications to statistics. In {\em
  Modern methods of multivariate statistics}; Hermann,  2014; Vol.~82, pp.
  135--154.

\bibitem[Ishi(2015)]{Is15}
Ishi, H.
\newblock Matrix realization of a homogeneous cone. In {\em Geometric science
  of information}; Springer, Cham,  2015; Vol. 9389, {\em Lecture Notes in
  Comput. Sci.}, pp. 248--256.

\end{thebibliography}

\end{document}